\documentclass[12pt,reqno]{amsart}

\usepackage{amsmath, epsfig, cite,color}
\usepackage{amssymb}
\usepackage{amsfonts}
\usepackage{latexsym}
\usepackage{graphicx}
\usepackage{algorithmic,algorithm}
\usepackage{amsthm,array,multirow,longtable}
\usepackage{tikz,caption}
\usetikzlibrary{decorations.pathreplacing}

\numberwithin{equation}{section}

\newcommand{\nc}{\newcommand}
\nc{\sha}{\mbox{\cyr X}}  %used to be \cyr
\font\cyr=wncyr10

\makeatletter \@addtoreset{equation}{section} \makeatother

\begin{document}
\allowdisplaybreaks
\title
{A bijection between two subfamilies of Motzkin paths}
\author[N.S.S. Gu]{Nancy S.S. Gu}
\address[N.S.S. Gu]{Center for Combinatorics, LPMC, Nankai
University, Tianjin 300071, P.R. China} \email{gu@nankai.edu.cn}
\author[H. Prodinger]{Helmut Prodinger}
\address[H. Prodinger]{Department of Mathematics, University of Stellenbosch,
7602 Stellenbosch, South Africa} \email{hproding@sun.ac.za}



\date{\today}

\begin{abstract} Two subfamilies of Motzkin paths, with the same numbers of up, down, horizontal steps were known to be equinumerous with ternary trees and related objects. We construct a bijection between these two families that does not use any auxiliary objects, like ternary trees.

\end{abstract}

\maketitle

\section{Introduction}

Motzkin paths are similar to Dyck paths, but allow also horizontal steps of unit length. In this note, we concentrate on two
subfamilies, where there are $n$ up steps ($u$), $n$ down steps ($d$), and $n$ horizontal steps ($h$). Both families are enumerated by
\begin{equation*}
1+ z^3+ 3z^6 +12z^9+ 55z^{12}+ 273z^{15}+\cdots,
\end{equation*}
and the coefficients also enumerate ternary trees and many other objects, see sequence A001764 in \cite{OEIS}.

The first family originates from Asinowski and Mansour \cite{AM}. They start from a Dyck path of length $2n$ and \emph{label} each maximal sequence of up steps by a Dyck path. If we say \emph{replace} instead of \emph{label} and use the steps $h$ and $u$ for the replaced sequence, we have a Motzkin path with $n$ horizontal steps. To clarify, we give a list of all 12 such paths of length 9:

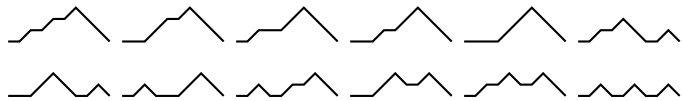
\begin{center}
	\begin{tikzpicture}[scale=0.15]
	
	%\draw[step=1.cm,black,dotted] (-0.0,-0.0) grid (9.0,3.0);

	\draw[thick] (0,0) to (1,0) to (2,1) to (3,1) to (4,2) to (5,2) to (6,3) to (7,2) to (8,1) to (9,0);
	
	\end{tikzpicture}
	\begin{tikzpicture}[scale=0.15]
	
	%\draw[step=1.cm,black,dotted] (-0.0,-0.0) grid (9.0,3.0);

	\draw[thick] (0,0) to (1,0) to (2,0) to (3,1) to (4,2) to (5,2) to (6,3) to (7,2) to (8,1) to (9,0);
	
	\end{tikzpicture}
	\begin{tikzpicture}[scale=0.15]
	
	%\draw[step=1.cm,black,dotted] (-0.0,-0.0) grid (9.0,3.0);

	\draw[thick] (0,0) to (1,0) to (2,1) to (3,1) to (4,1) to (5,2) to (6,3) to (7,2) to (8,1) to (9,0);
	
	\end{tikzpicture}
	\begin{tikzpicture}[scale=0.15]
	
	%\draw[step=1.cm,black,dotted] (-0.0,-0.0) grid (9.0,3.0);

	\draw[thick] (0,0) to (1,0) to (2,0) to (3,1) to (4,1) to (5,2) to (6,3) to (7,2) to (8,1) to (9,0);
	
	\end{tikzpicture}
	\begin{tikzpicture}[scale=0.15]
	
	%\draw[step=1.cm,black,dotted] (-0.0,-0.0) grid (9.0,3.0);

	\draw[thick] (0,0) to (1,0) to (2,0) to (3,0) to (4,1) to (5,2) to (6,3) to (7,2) to (8,1) to (9,0);
	
	\end{tikzpicture}
	\begin{tikzpicture}[scale=0.15]
	
	%\draw[step=1.cm,black,dotted] (-0.0,-0.0) grid (9.0,3.0);

	\draw[thick] (0,0) to (1,0) to (2,1) to (3,1) to (4,2) to (5,1) to (6,0) to (7,0) to (8,1) to (9,0);
	
	\end{tikzpicture}
	
	\begin{tikzpicture}[scale=0.15]
	
	%\draw[step=1.cm,black,dotted] (-0.0,-0.0) grid (9.0,3.0);

	\draw[thick] (0,0) to (1,0) to (2,0) to (3,1) to (4,2) to (5,1) to (6,0) to (7,0) to (8,1) to (9,0);
	
	\end{tikzpicture}
	\begin{tikzpicture}[scale=0.15]
	
	%\draw[step=1.cm,black,dotted] (-0.0,-0.0) grid (9.0,3.0);

	\draw[thick] (0,0) to (1,0) to (2,1) to (3,0) to (4,0) to (5,0) to (6,1) to (7,2) to (8,1) to (9,0);
	
	\end{tikzpicture}
	\begin{tikzpicture}[scale=0.15]
	
	%\draw[step=1.cm,black,dotted] (-0.0,-0.0) grid (9.0,3.0);

	\draw[thick] (0,0) to (1,0) to (2,1) to (3,0) to (4,0) to (5,1) to (6,1) to (7,2) to (8,1) to (9,0);
	
	\end{tikzpicture}
	\begin{tikzpicture}[scale=0.15]
	
	%\draw[step=1.cm,black,dotted] (-0.0,-0.0) grid (9.0,3.0);

	\draw[thick] (0,0) to (1,0) to (2,0) to (3,1) to (4,2) to (5,1) to (6,1) to (7,2) to (8,1) to (9,0);
	
	\end{tikzpicture}
	\begin{tikzpicture}[scale=0.15]
	
	%\draw[step=1.cm,black,dotted] (-0.0,-0.0) grid (9.0,3.0);

	\draw[thick] (0,0) to (1,0) to (2,1) to (3,1) to (4,2) to (5,1) to (6,1) to (7,2) to (8,1) to (9,0);
	
	\end{tikzpicture}
	\begin{tikzpicture}[scale=0.15]
	
	%\draw[step=1.cm,black,dotted] (-0.0,-0.0) grid (9.0,3.0);

	\draw[thick] (0,0) to (1,0) to (2,1) to (3,0) to (4,0) to (5,1) to (6,0) to (7,0) to (8,1) to (9,0);
	
	\end{tikzpicture}
\end{center}

The second family was introduced to model \emph{frog hops} from a question in a student's olympiad \cite{bulgaria}:
They were called S-Motzkin paths, have the same number of $u$, $d$, $h$, and when deleting the down steps, the sequence must look like $huhuhu\dots hu$. Here is again a list of all 12 objects of length 9:

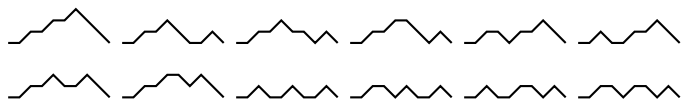
\begin{center}
	\begin{tikzpicture}[scale=0.15]
	
	%\draw[step=1.cm,black,dotted] (-0.0,-0.0) grid (9.0,3.0);

	\draw[thick] (0,0) to (1,0) to (2,1) to (3,1) to (4,2) to (5,2) to (6,3) to (7,2) to (8,1) to (9,0);
	
	\end{tikzpicture}
	\begin{tikzpicture}[scale=0.15]
	
	%\draw[step=1.cm,black,dotted] (-0.0,-0.0) grid (9.0,3.0);

	\draw[thick] (0,0) to (1,0) to (2,1) to (3,1) to (4,2) to (5,1) to (6,0) to (7,0) to (8,1) to (9,0);
	
	\end{tikzpicture}
	\begin{tikzpicture}[scale=0.15]
	
	%\draw[step=1.cm,black,dotted] (-0.0,-0.0) grid (9.0,3.0);

	\draw[thick] (0,0) to (1,0) to (2,1) to (3,1) to (4,2) to (5,1) to (6,1) to (7,0) to (8,1) to (9,0);
	
	\end{tikzpicture}
	\begin{tikzpicture}[scale=0.15]
	
	%\draw[step=1.cm,black,dotted] (-0.0,-0.0) grid (9.0,3.0);

	\draw[thick] (0,0) to (1,0) to (2,1) to (3,1) to (4,2) to (5,2) to (6,1) to (7,0) to (8,1) to (9,0);
	
	\end{tikzpicture}
	\begin{tikzpicture}[scale=0.15]
	
	%\draw[step=1.cm,black,dotted] (-0.0,-0.0) grid (9.0,3.0);

	\draw[thick] (0,0) to (1,0) to (2,1) to (3,1) to (4,0) to (5,1) to (6,1) to (7,2) to (8,1) to (9,0);
	
	\end{tikzpicture}
	\begin{tikzpicture}[scale=0.15]
	
	%\draw[step=1.cm,black,dotted] (-0.0,-0.0) grid (9.0,3.0);

	\draw[thick] (0,0) to (1,0) to (2,1) to (3,0) to (4,0) to (5,1) to (6,1) to (7,2) to (8,1) to (9,0);
	
	\end{tikzpicture}
	
	\begin{tikzpicture}[scale=0.15]
	
	%\draw[step=1.cm,black,dotted] (-0.0,-0.0) grid (9.0,3.0);

	\draw[thick] (0,0) to (1,0) to (2,1) to (3,1) to (4,2) to (5,1) to (6,1) to (7,2) to (8,1) to (9,0);
	
	\end{tikzpicture}
	\begin{tikzpicture}[scale=0.15]
	
	%\draw[step=1.cm,black,dotted] (-0.0,-0.0) grid (9.0,3.0);

	\draw[thick] (0,0) to (1,0) to (2,1) to (3,1) to (4,2) to (5,2) to (6,1) to (7,2) to (8,1) to (9,0);
	
	\end{tikzpicture}
	\begin{tikzpicture}[scale=0.15]
	
	%\draw[step=1.cm,black,dotted] (-0.0,-0.0) grid (9.0,3.0);

	\draw[thick] (0,0) to (1,0) to (2,1) to (3,0) to (4,0) to (5,1) to (6,0) to (7,0) to (8,1) to (9,0);
	
	\end{tikzpicture}
	\begin{tikzpicture}[scale=0.15]
	
	%\draw[step=1.cm,black,dotted] (-0.0,-0.0) grid (9.0,3.0);

	\draw[thick] (0,0) to (1,0) to (2,1) to (3,1) to (4,0) to (5,1) to (6,0) to (7,0) to (8,1) to (9,0);
	
	\end{tikzpicture}
	\begin{tikzpicture}[scale=0.15]
	
	%\draw[step=1.cm,black,dotted] (-0.0,-0.0) grid (9.0,3.0);

	\draw[thick] (0,0) to (1,0) to (2,1) to (3,0) to (4,0) to (5,1) to (6,1) to (7,0) to (8,1) to (9,0);
	
	\end{tikzpicture}
	\begin{tikzpicture}[scale=0.15]
	
	%\draw[step=1.cm,black,dotted] (-0.0,-0.0) grid (9.0,3.0);

	\draw[thick] (0,0) to (1,0) to (2,1) to (3,1) to (4,0) to (5,1) to (6,1) to (7,0) to (8,1) to (9,0);
	
	\end{tikzpicture}
\end{center}

These paths and some of their properties were investigated in the recent paper \cite{PSW}.

The goal of the present note is to describe a bijection between the two families, which operates strictly on the paths themselves, without involving any other objects that are equinumerous.

Instead of drawing pictures, we use the more economical description with the letters $u,d,h$, and we
always consider the paths from left to right.

\section{From paths of the Asinowski/Mansour type to S-Motzkin paths}

From paths given by Asinowski and Mansour, let $A_0 D_1 A_1 D_2 \ldots A_{t-1}D_t$ denote a path where $A_i$ ($i=0,1,\dots,t-1$) denotes a sub-path consisting of horizontal and up steps, and $D_i$ ($i=1,2,\dots,t$) denotes a sub-path consisting of consecutive down steps. Let $d_i$ ($i=1,2,\dots,t$) denote the number of steps of $D_i$. For completeness, we mention that all $A_i$ and $D_i$ are non-empty.

First, we replace horizontal steps and up steps by up steps and down steps, respectively, in $A_i$ to get a Dyck path $\overline{A_i}$. Then for each $\overline{A_i}$, we insert horizontal steps into it to form an S-Motzkin path $A_i'$ such that each horizontal step is in the first position that can be inserted from left to right. Observe that in $A_i'$, there are no two horizontal steps on the same height such that all the steps between these two steps are above that height, and each up step except for the last one must be followed (eventually) by a horizontal step.

These canonically obtained S-Motzkin paths are the \emph{building stones} of the final object, and the numbers $d_1,d_2,\dots,d_{t-1}$ are used
to tell us how to \emph{glue} them together. Notice also that the number $d_t$ will not be used for the construction.

Let $u_i$ denote the number of up steps of $A_i'$. Notice that
\begin{align}
\sum_{k=1}^i d_k & \leq \sum_{k=0}^{i-1} u_k \quad  (1 \leq i \leq t-1),\label{t-1}\\
\sum_{k=1}^t d_k &= \sum_{k=0}^{t-1} u_k.\label{t}
\end{align}

First, we draw $A_0'$. Then according to $d_1,d_2,\dots,d_{t-1}$, we insert $A_1,A_2,\dots,A_{t}$ in turn. If $d_1 < u_0$, then from the beginning of $A_0'$, find the $d_1$-th up step and insert $A_1'$ behind it. If $d_1=u_0$, then insert $A_1'$ behind the last step of $A_0'$.

Assume that we have inserted $A_{i-1}$ ($1\leq i \leq t-1$). Then from the beginning of $A_{i-1}'$, find the $d_i$-th up step. Notice that if $d_i <u_{i-1}$, then the $d_i$-th up step belongs to $A_{i-1}'$ and is not the last up step of $A_{i-1}'$. We insert $A_i'$ behind this $d_i$-th up step. If $d_i=u_{i-1}$, then insert $A_i'$ behind the end of $A_{i-1}'$. If $d_i>u_{i-1}$, then the $d_i$-th up step may belong to $A_0'$, $A_1'$, $\dots$, or $A_{i-2}'$. This time, if this step is not any last up step of $A_0'$, $A_1'$, $\dots$, or $A_{i-2}'$, we insert $A_i'$ behind it. Otherwise, if this step is the last up step of $A_s'$ ($s\in \{0,1,\dots,i-2\}$), then we put $A_i'$ behind $A_s'$. Note that \eqref{t-1} and \eqref{t} ensure that we can always find the aimed up step.

Finally, we get an S-Motzkin path until we have inserted all $A_i'$ ($i=0,1,\dots,t-1$).

\section{From S-Motzkin paths to paths of the Asinowski/Mansour type}

For an S-Motzkin path $P$, from left to right, we check the horizontal steps in turn. For a given horizontal step $h_i$, look along the path from $h_i$. If there is a closest horizontal step on the same height and all the steps between these two steps denoted by $P_i$ are not below this height, then $h_iP_i$ is an S-Motzkin path. We call $h_i$ a paired horizontal step. For the first horizontal step of the S-Motzkin path $P$, no matter whether it is paired or not, we always call it a paired horizontal step denoted by $h_0$. If there are no other horizontal steps on height $0$, then $P=h_0P_0$. Let $h_0,h_1,\dots,h_{t-1}$ denote these paired horizontal steps of $P$.

If there are no paired horizontal steps in $P_i$, then $h_iP_i$ is an S-Motzkin sub-path with all the horizontal steps located in their first positions. We assume $h_iP_i$ as $A_i'$. Back to the S-Motzkin path $P$, if the up step just before $h_i$ is the $d_i$-th up step when we are counting from the beginning of $A_{i-1}'$, then we record $d_i$. Then in the S-Motzkin path $P$, we find the next paired horizontal step behind $h_i$, and repeat the process.

If there is a paired horizontal step $h_{i+1}$ in $P_i$, then we use $h_{i+1}$ to find $A_{i+1}'$. Note that
now $A_i'$ contains the steps of $h_iP_i$ without those steps in $A_{i+1}'$.

In view of all the paired horizontal steps, and using the above method, we can identify all $A_i'$ ($i=0,1,\dots,t-1$) and $d_i$ ($i=1,2,\dots,t-1$).

Deleting all the horizontal steps in $A_i'$ ($0 \leq i \leq t-1$), and then replacing up and down steps by horizontal and up steps, respectively, we obtain $A_i$. Let $D_i$ ($1 \leq i \leq t-1$) be the sub-path consisting of $d_i$ consecutive down steps. Then we draw a path $A_0 D_1 A_1 D_2 A_2 \dots D_{t-1} A_{t-1}$, and add enough consecutive down steps after it to form a Motzkin path given by Asinowski and Mansour.

\section{A detailed example}

We consider  a path given by Asinowski and Mansour
$$hhhuuu\mathbf{d}hhuuhu\mathbf{d}hu\mathbf{dd}hhuu\mathbf{ddddd}.$$
Then we can divide the paths into four parts $A_0$, $A_1$, $A_2$ and $A_3$ by using the consecutive down steps which are in bold. Set
\begin{align*}
\begin{array}{llll}
A_0=hhhuuu,& \qquad A_1=hhuuhu,& \qquad A_2=hu,& \qquad A_3=hhuu,\\
d_1=1,& \qquad d_2=1, &\qquad d_3=2, &\qquad d_4=5.
\end{array}
\end{align*}
In fact, we do not need $d_4$ in the bijection.

First, replacing horizontal steps and up steps by up steps and down steps, respectively, in $A_i$ ($i=0,1,2,3$), we have four Dyck paths
\begin{align*}
\overline{A_0}=uuuddd,\qquad \overline{A_1}=uuddud,\qquad \overline{A_2}=ud,\qquad \overline{A_3}=uudd.
\end{align*}
Then inserting horizontal steps in $\overline{A}_i$ ($i=0,1,2,3$) in the first possible positions, we get
\begin{align*}
A_0'=huhuhuddd,\qquad A_1'=huhuhddud, \qquad A_2'=hud, \qquad A_3'=huhudd.
\end{align*}
Since $d_1=1$, we find the first up step in $A_0'$. Then inserting $A_1'$ behind this step, we obtain
$$hu\mathbf{huhuhddud}huhuddd.$$
For the above path, since $d_2=1$, we find the first up step from the beginning of $A_1'$. That is to say, deleting the first two steps in the above path, we get a sub-path. For this sub-path, we find the $d_2$-th up step, and then insert $A_2'$. Let $u_1$ be the number of up steps of $A_1'$. Note that if $d_2<u_1$, then $A_2'$ is inserted in $A_1'$. If $d_2=u_1$, then put $A_2'$ just behind $A_1'$. If $d_2>u_1$, then $A_2$ is inserted between two steps of $A_0'$. Back to this example, we have $d_2<u_1$. So
$$huhu\mathbf{hud}huhddudhuhuddd.$$
For the above path, since $d_3=2$, we find the second up step from the beginning of $A_2'$. That is to say, deleting the first four steps in the above path yields a sub-path. In this sub-path, we find the second up step. In this example, $d_3>u_2$, and the second up step is not the last up step of $A_1'$ or $A_0'$. So we insert $A_3'$ behind the second up step directly. we have
$$huhuhudhu\mathbf{huhudd}hddudhuhuddd$$
which is an S-Motzkin path.

Inversely, for the S-Motzkin path
$$huhuhudhuhuhuddhddudhuhuddd,$$
we can find the paired horizontal steps as follows:
$$\mathbf{h}u\mathbf{h}u\mathbf{h}udhu\mathbf{h}uhuddhddudhuhuddd.$$
The first horizontal step must be a paired horizontal step although there is no other horizontal step on level $0$. We mark it as 
$$(\bar{h}uhuhudhuhuhuddhddudhuhuddd).$$
Here we use a pair of parentheses to identify $A_0'$. 

Then we find the next pair of horizontal steps
$$(\bar{h}u(\check{h}uhudhuhuhuddhddud)\check{h}uhuddd).$$
So we have $d_1=1$, and we use a pair of parentheses to identify $A_1'$ from $A_0'$. 

Next, we have
$$(\bar{h}u(\check{h}u(\hat{h}ud)\hat{h}uhuhuddhddud)\check{h}uhuddd)$$
and $d_2=1$. We add a pair of parentheses to identify $A_2'$. 

Finally, we obtain
$$(\bar{h}u(\check{h}u(\hat{h}ud)\hat{h}u(\dot{h}uhudd)\dot{h}ddud)\check{h}uhuddd)$$
and $d_3=2$ where we use a pair of parentheses to identify $A_3'$.

Therefore, we have
\begin{align*}
\begin{array}{llll}
A_0'=huhuhuddd, &\qquad A_1'=huhuhddud, &\qquad A_2'=hud, &\qquad A_3'=huhudd,\\*
d_1=1,&\qquad  d_2=1,&\qquad  d_3=2.&
\end{array}
\end{align*}

Deleting all the horizontal steps in $A_0'$, $A_1'$, $A_2'$ and $A_3'$, and then replacing up and down steps by horizontal and up steps, respectively, we obtain
\begin{align*}
A_0&=hhhuuu,\qquad A_1=hhuuhu, \qquad A_2=hu,\qquad A_3=hhuu.
\end{align*}
Combining $A_i$ ($i=0,1,2,3$) and $d_i$ ($i=1,2,3$), we have the following path
$$hhhuuu\mathbf{d}hhuuhu\mathbf{d}hu\mathbf{dd}hhuu.$$
Finally, we add enough down steps at the end of the above path to form a Motzkin path given by Asinowski and Mansour:
$$hhhuuu\mathbf{d}hhuuhu\mathbf{d}hu\mathbf{dd}hhuu\mathbf{ddddd}.$$

\section{The objects of length 9 matched}

For the reader's convenience, we provide the correspondence of 12 objects. The number 12 is very convenient; it is not too small and not too large. In \cite{PSW}, there were also many explicit lists with 12 objects each.
\begin{center}
	
	\begin{table}[h]
	\begin{tabular}{r | l | l || l}
		
		No. & AM-paths & AM-paths dec.& S-Motzkin paths\\
		
		\hline
		
		1&$huhuhuddd$&$ududud$\ 3&$huhduhdud$\\
		
		2&$hhuuhuddd$&$uuddud$\ 3&$huhuhddud$\\
		
		3&$huhhuuddd$&$uduudd$\ 3&$huhduhudd$\\
		
		4&$hhuhuuddd$&$uududd$\ 3&$huhuhdudd$\\
		
		5&$hhhuuuddd$&$uuuddd$\ 3&$huhuhuddd$\\
		
		6&$huhuddhud$&$udud$\ 2\ $ud$\ 1&$huhdudhud$\\
		
		7&$hhuuddhud$&$uudd$\ 2\ $ud$\ 1&$huhuddhud$\\
		
		8&$hudhhuudd$&$ud$\ 1\ $uudd$\ 2&$hudhuhudd$\\
		
		9&$hudhuhudd$&$ud$\ 1\ $udud$\ 2&$hudhuhdud$\\
		
		10&$hhuudhudd$&$uudd$\ 1\ $ud$\ 2&$huhudhudd$\\
		
		11&$huhudhudd$&$udud$\ 1\ $ud$\ 2&$huhudhdud$\\
		
		12&$hudhudhud$&$ud$\ 1\ $ud$\ 1\ $ud$\ 1&$hudhudhud$\\

	\end{tabular}

	\caption{Paths from Asinowski and Mansour, also decomposed, and the corrresponding S-Motzkin paths.}
\end{table}	
\end{center}

\noindent {\bf Acknowledgements:} This work was supported by the National Natural Science Foundation of China and the
Fundamental Research Funds for the Central Universities (Nankai University).

\end{document}